\documentclass{amsart}
\usepackage{graphicx}
\usepackage{amssymb}
\usepackage{amsfonts}
\swapnumbers
\newtheorem{theorem}{Theorem}[section]

\newtheorem{corollary}[theorem]{Corollary}
\newtheorem{proposition}[theorem]{Proposition}
\theoremstyle{definition}

\newtheorem{assumption}[theorem]{Assumption}

\numberwithin{equation}{section}
 \theoremstyle{plain}    
 
 \numberwithin{equation}{section} 
 \numberwithin{figure}{section} 
 \theoremstyle{plain}    
 \theoremstyle{plain}    
 \theoremstyle{remark}    
 \newtheorem*{acknowledgement*}{Acknowledgement} 

\newcommand{\cA}{{\mathcal A}}

\newcommand{\cF}{{\mathcal F}}
\newcommand{\cG}{{\mathcal G}}
\newcommand{\cH}{{\mathcal H}}

\newcommand{\cM}{{\mathcal M}}

\newcommand{\te}{{\theta}}

\newcommand{\Om}{{\Omega}}
\newcommand{\om}{{\omega}}
\newcommand{\ve}{{\varepsilon}}
\newcommand{\del}{{\delta}}

\newcommand{\gam}{{\gamma}}

\newcommand{\vp}{{\varpi}}
\newcommand{\io}{{\iota}}
\newcommand{\up}{{\upsilon}}

\newcommand{\sig}{{\sigma}}
\newcommand{\al}{{\alpha}}
\newcommand{\be}{{\beta}}
\newcommand{\ka}{{\kappa}}
\newcommand{\la}{{\lambda}}

\newcommand{\bbR}{{\mathbb R}}

\newcommand{\bbI}{{\mathbb I}}

\begin{document}
\title[]{A strong invariance principle\\
 for nonconventional sums}%
 \vskip 0.1cm 
 \author{ Yuri Kifer\\
 \vskip 0.1cm
Institute of Mathematics\\
The Hebrew University of Jerusalem}%
\email{kifer\@@math.huji.ac.il}
\address{Institute of Mathematics, Hebrew University, Jerusalem 91904,\linebreak
 Israel}

\thanks{ }
\dedicatory{  }
 \date{\today}
\begin{abstract}\noindent
In \cite{KV} we obtained a functional central limit theorem  (known also
as a weak invariance principle) for sums of the form 
$\sum_{n=1}^{[Nt]}F\big(X(n),X(2n),...,X(kn),
X(q_{k+1}(n)),X(q_{k+2}(n)),\cdots, X(q_\ell(n))\big)$
(normalized by $1/\sqrt N$) where $X(n),n\geq 0$ is a sufficiently fast 
mixing vector process with some moment conditions and stationarity properties,
$F$ is a continuous function with polinomial growth and certain regularity 
properties and $q_i,i>m$ are positive functions taking on integer values on
integers with some growth conditions which are satisfied, for instance, when
$q_i$'s are polynomials of growing degrees. This paper deals with strong
invariance principles (known also as strong approximation theorems) for such
sums which provide their uniform in time almost sure approximation by processes
built out of Brownian motions with error terms growing
slower than $\sqrt N$. This yields, in particular, an invariance principle 
in the law of iterated algorithm for the above sums.
Among motivations for such results are their applications to multiple
recurrence for stochastic processes and dynamical systems as well, as to some
questions in metric number theory and they can be considered as a natural 
follow up of a series of papers dealing with nonconventional ergodic averages.

\end{abstract}
\maketitle
\markboth{Yu.Kifer}{Strong invariance principles} 
\renewcommand{\theequation}{\arabic{section}.\arabic{equation}}
\pagenumbering{arabic}

\section{Introduction}\label{sec1}\setcounter{equation}{0}

Nonconventional ergodic theorems attracted substantial attention 
in ergodic theory (see, for instance, \cite{Be} and \cite{Fu}). From 
a probabilistic point of view ergodic theorems are laws of large numbers
for stationary processes and once they are established it is natural to 
study deviations from the average. The most celebrated result of this kind
is the central limit theorem. In \cite{KV} we obtained a functional limit
theorem for expressions of the form
\begin{equation}\label{1.1}
\xi_N(t)=1/\sqrt {N}\sum_{1\leq n\leq Nt}\big(F\big(X(q_1(n)),..., X(q_\ell(n))
\big)-\bar F\big)
\end{equation}
and for the corresponding continuous time expressions of the form
\begin{equation}\label{1.2}
\xi_N(t)=1/\sqrt {N}\int_0^{Nt}\big(F\big( X(q_1(t)),..., X(q_\ell(t))\big)
-\bar F\big)dt
\end{equation} 
where $X(n),n\geq 0$'s is a sufficiently fast mixing 
vector valued process with some moment conditions and stationarity properties,
$F$ is a continuous function with polinomial growth and certain regularity
properties, 
$\bar F=\int Fd(\mu\times\cdots\times\mu)$, $\mu$ is the distribution
of $X(0)$, $q_j(t)=jt,\, j\leq k$ and $q_j,j>k$ are positive functions 
taking on integer values on integers in the discrete time case
with some growth conditions which are satisfied, for instance, when
$q_i$'s are polynomials of growing degrees though we actully need much less. 
A substantially more restricted central limit theorem for expressions of this
 sort was obtained in \cite{Ki}.
 
 Functional central limit theorems are called nowadays also weak invariance
 principles while for more than 40 years now (since probably Strassen's
 work \cite{St}) probabilists were interested also in strong invariance
 principles called also strong approximation theorems. The latter provides
 almost sure or in average approximation of a sum of $N$ random variables
 by a Brownian motion or, more generally, by a Gaussian process with an
 error term growing slower than $\sqrt N$ which yields as a result both 
 the central limit theorem and the law of iterated logarithm, as well as
other limiting results which are clear or easy to prove for Gaussian
processes. 

We will show in this paper that the sums $\Xi(Nt)=\sqrt {N}\xi_N(t)$ 
appearing in (\ref{1.1}) can be represented as $\sum_{1\leq i\leq\ell}
\Xi_i(Nt)$ where each $\Xi_i(Nt)$ can be approximated with an error term
of order $N^{\frac 12-\al},\al>0$ by a process $\sig_iB_i(t)$
where $\sig_i\geq 0$ is a constant and $B_i$ is a Brownian motion.
This result yields also a law of iterated logarithm type result saying that
with probability one all limit points as $N\to\infty$ of the sequence 
$\xi_N(t)(\log\log N)^{-1/2},\,t\in[0,1]$ of paths belong to a compact set.

Our methods employ the martingale approximation machinery from \cite{KV},
enhanced so that to obtain appropriate error estimates, together with the
technique from \cite{PS} which involves partition into blocks and Skorokhod
embedding of martingales into a Brownian motion (the latter was first used
for similar purposes in \cite{St}). Observe that the summands in (\ref{1.1})
 depend strongly on the future and martingale methods start working only
 after we force "the future to become present". By this reason the role
of martingales in our nonconventional framework was not selfevident at the
beginning but their effective use initiated in \cite{KV} opened a wide vista
for proving various limit theorems in this setup. 
It was shown in \cite{KV} that $\xi_N$ converges weakly to a Gaussian
process and it would be interesting to obtain a strong approximation of
$\sqrt {N}\xi_N(t)$ by such Gaussian process but this would require to
deal with multi dimensional approximations where the Skorokhod embedding
we rely on does not work. Observe that since the 1960ies 
several other methods were developed to provide approximation of sums of 
random variables by a Brownian motion. Among them is the quantile method (see,
for instance, \cite{KMT}) which provides essentially optimal approximation 
but works only for independent random variables and by this reason does
not seem applicable to our setup. Another method developed by Stein
(see its recent account in \cite{Ch}) also yields nearly optimal error
estimates but it is not yet clear whether it can be adapted to our situation.
The advantages of yet another method based on estimates of conditional
characteristic functions (see, for instance, \cite{BP}) lie in its
applicability to the multidimensional situation where, for instance, the
Skorokhod embedding does not work well, but complications in the use of
characteristic functions exhibited in \cite{Ki} make applicability of
this method in our setup doubtful. 

As in \cite{KV} our results hold true when, for instance, $X(n)=T^nf$ where
$f=(f_1,...,f_\wp)$, $T$ is a mixing
 subshift of finite type, a hyperbolic 
diffeomorphism or an expanding transformation taken with a Gibbs invariant 
measure (see for instance, \cite{Bo}) and some other
 dynamical systems, as well, as in
 the case when $X(n)=f(\xi_n),\, 
 f=(f_1,...,f_\wp)$ where $\xi_n$ is a Markov chain
 satisfying the
 Doeblin condition (see \cite{IL}) considered as a stationary process
 with
  respect to its invariant measure. The main known application of the above
  type results is to multiple recurrence when we employ our limit theorems
  for the random variable which counts returns of the stochastic process
  under consideration to given sets. In this case the function $F$ above
  is a product of some coordinate functions in which we plug in 
  corresponding $X(q_i(n))=\bbI_{A_j}(\eta(q_i(n)))$ where $\eta(m)$ is
  either $T^mx$ in the dynamical systems case or $\xi_m$ in the Markov
  chain case and $\bbI_A$ is the indicator of a set $A$. This yields also
  applications to metric number theory providing limit theorems, for instance,
  for the number $M_N(x)$ of times first $N$ digits in the $m$-base or
  continued fraction expansion of $x$ belong to a chosen subset of digits.
  As it is well known the former expansions can be obtained via the 
  multiplication by $m$ (expanding) transformation while the latter via
  the Gauss map of the interval and both dynamical systems are 
  exponentially fast $\psi$ mixing with respect to their invariant Lebesgue
  or Gauss measure, respectively (see, for instance, \cite{He}).

  \section{Preliminaries and main results}\label{sec2}\setcounter{equation}{0}

Our setup consists of a  $\wp$-dimensional stochastic process
$\{X(n),  n=0,1,...\}$ on a probability space $(\Om,\cF,P)$ and of a family
of $\sig$-algebras \hbox{$\cF_{kl}\subset\cF,\, -\infty\leq k\leq
l\leq\infty$} such that $\cF_{kl}\subset\cF_{k'l'}$ if $k'\leq k$ and
$l'\geq l$. The dependence between two sub $\sig$-algebras $\cG,\cH\subset\cF$
is measured often via the quantities
\begin{equation}\label{2.1}
\varpi_{q,p}(\cG,\cH)=\sup\{\| E\big [g|\cG\big]-E[g]\|_p:\, g\,\,\mbox{is}\,\,
\cH-\mbox{measurable and}\,\,\| g\|_q\leq 1\},
\end{equation}
where the supremum is taken over real functions and $\|\cdot\|_r$ is the
$L^r(\Om,\cF,P)$-norm. Then more familiar $\al,\rho,\phi$ and $\psi$-mixing 
(dependence) coefficients can be expressed in the form (see \cite{Bra}, Ch. 4 ),
\begin{eqnarray*}
&\al(\cG,\cH)=\frac 14\varpi_{\infty,1}(\cG,\cH),\,\,\rho(\cG,\cH)=\varpi_{2,2}
(\cG,\cH)\\
&\phi(\cG,\cH)=\frac 12\varpi_{\infty,\infty}(\cG,\cH)\,\,\mbox{and}\,\,
\psi(\cG,\cH)=\varpi_{1,\infty}(\cG,\cH).
\end{eqnarray*}
The relevant quantities in our setup are
\begin{equation}\label{2.2}
\varpi_{q,p}(n)=\sup_{k\geq 0}\varpi_{q,p}(\cF_{-\infty,k},\cF_{k+n,\infty})
\end{equation}
and accordingly
\[
\al(n)=\frac{1}{4}\varpi_{\infty,1}(n),\,\rho(n)=\varpi_{2,2}(n),\,
\phi(n)=\frac 12\varpi_{\infty,\infty}(n)\,\,\mbox{and}\,\,\ \psi(n)=
\varpi_{1,\infty}(n).
\]
Our assumptions will require certain speed of decay as $n\to\infty$ of both
the mixing rates $\varpi_{q,p}(n)$ and the approximation rates defined by
\begin{equation}\label{2.3}
\beta_p (n)=\sup_{m\geq 0}\|X(m)-E\big (X(m)|\cF_{m-n,m+n}\big)\|_p.
\end{equation}
In what follows we can always extend the definitions of $\cF_{kl}$
given only for $k,l\geq 0$ to negative $k$ by defining  $\cF_{kl}=\cF_{0l}$ 
for $k<0$ and $l\geq 0$. Furthermore, we do not require stationarity of the
process $X(n), n\geq 0$ assuming only that the distribution of $X(n)$ does not
depend on $n$ and the joint distribution of $\{X(n), X(n')\}$  depends only on
$n-n'$ which we write for further references by 
\begin{equation}\label{2.4}
X(n)\stackrel {d}{\sim}\mu\,\,\mbox{and}\,\,
(X(n),X(n'))\stackrel {d}{\sim}\mu_{n-n'}\,\,\mbox{for all}\,\,n,n'
\end{equation}
where $Y\stackrel {d}{\sim}Z$ means that $Y$ and $Z$ have the same distribution.

Next,let $F= F(x_1,...,x_\ell),\, x_j\in\bbR^{\wp}$ be a function on 
$\bbR^{\wp\ell}$ such that for some $\iota,K>0,\ka\in (0,1]$ and all 
$x_i,y_i\in\bbR^{\wp}, i=1,...,\ell$, 
\begin{equation}\label{2.5}
|F(x_1,...,x_\ell)-F(y_1,...,y_\ell)|\leq K\big(1+\sum^\ell_{j=1}|x_j|^\iota+
\sum^\ell_{j=1} |y_j|^\iota \big)\sum^\ell_{j=1}|x_j-y_j|^\ka
\end{equation}
and 
\begin{equation}\label{2.6} 
|F(x_1,...,x_\ell)|\leq K\big( 1+\sum^\ell_{j=1}|x_j|^{\iota} \big).
\end{equation}
The above assumptions are motivated by the desire to include, for instance,
functions $F$ polinomially dependent on their arguments.
To simplify formulas we assume a centering condition
\begin{equation}\label{2.7}
\bar F=\int F(x_1,...,x_\ell)\,d\mu(x_1)\cdots d\mu(x_\ell)=0
\end{equation}
which is not really a restriction since we always can replace $F$ by 
$F-\bar F$.

Our setup includes also a sequence of increasing functions $q_1(n)< q_2(n) 
<\cdots < q_\ell(n)$ taking on integer values on integers and such that the 
first $k$ of them are $q_j(n)=jn,\, j\leq k$ whereas the remaining ones grow
faster in $n$. We assume that for $k+1\le i \le \ell$,
\begin{equation}\label{2.8}
q_i(n+1)-q_i(n)\geq n^\del
\end{equation}
for some $\del>0$ and all $n\geq 2$ while for $i\ge k$ and any $\epsilon > 0 $,
\begin{equation}\label{2.9}
\liminf_{n\to\infty}(q_{i+1}(\epsilon\, n)-q_i(n))>0
\end{equation}
which is equivalent in view of (\ref{2.8}) to
\begin{equation}\label{2.10}
\liminf_{n\to\infty}(q_{i+1}(\epsilon\, n)-q_i(n))=\infty.
\end{equation}

In order to give a detailed statement of our main result as well as for its
proof it will be essential to represent the function $F= F(x_1,x_2,\ldots,
x_\ell)$ in the form 
\begin{equation}\label{2.12}
F=F_1(x_1)+\cdots+F_\ell(x_1, x_2,\ldots, x_\ell)
\end{equation}
where for $i<\ell$,
\begin{eqnarray}\label{2.13}
&F_i(x_1,\ldots, x_i)=\int F(x_1,x_2, \ldots, x_\ell)\ d\mu (x_{i+1})\cdots 
d\mu(x_\ell)\\
&\quad -\int F(x_1,x_2, \ldots, x_\ell) \,d\mu (x_i)\cdots d\mu(x_\ell)\nonumber
\end{eqnarray}
and
\[
F_\ell(x_1,x_2, \ldots, x_\ell)=F(x_1,x_2, \ldots, x_\ell) -\int F(x_1,x_2, 
\ldots, x_\ell)\, d\mu(x_\ell)
\]
which ensures, in particular, that
\begin{equation}\label{2.14}
\int F_i(x_1, x_2,\ldots,x_{i-1}, x_i)\,d\mu(x_i)\equiv 0 \quad\forall 
\quad x_1, x_2,\ldots, x_{i-1}.
\end{equation}
These enable us to write
\begin{equation}\label{2.15}
\Xi(t)=\sum_{i=1}^\ell \Xi_i(t)
\end{equation}
where for $1\leq i\leq\ell $,
\begin{equation}\label{def2.17}
\Xi_i(t)=\sum_{1\leq n\leq t} F_i(X(q_1(n)),\ldots, X(q_i(n))).
\end{equation}
The decomposition of $\Xi(t)$ above is different from \cite{KV} since we
work here with each $\Xi_i(t)$ separately remaining all the time within
a one dimensional framework and do not care about multi dimensional 
covariances.

For each $\theta>0$ set
\begin{equation}\label{2.18}
\gamma_\theta^\theta = \|X\|_\theta^\theta= E|X(n)|^\theta  =
\int \|x\|^\theta d\mu .
\end{equation}
Our main result relies on
\begin{assumption}\label{ass2.1} With $d=(\ell-1)\wp$ there exist $p,q\ge 1$
 and $\delta,m >0$ with   $ \delta<\ka-\frac dp$ satisfying 
\begin{equation}\label{2.19}
\sum_{n=0}^\infty n^\del\varpi_{q,p}(n)<\infty,
\end{equation}
\begin{equation}\label{2.20}
\sum_{r=0}^\infty (r\beta_q(r))^\delta<\infty,
\end{equation}
\begin{equation}\label{2.21}
\gamma_{m}<\infty\,\,,\gam_{2q(\io+2)}<\infty\,\,\mbox{with}\,\,
 \frac{1}{2+\del}\ge \frac{1}{p}+\frac{\iota+2}{m}+\frac{\delta}{q}.
 \end{equation}
 \end{assumption}
Following \cite{PS} we will write $Z(t)\ll a(t)$ a.s. for a family of random
variables $Z(t),t\geq 0$ and a positive function $a(t),t\geq 0$ if
$\lim\sup_{t\to\infty}|Z(t)/a(t)|<\infty$ almost surely (a.s.)

\begin{theorem}\label{thm2.2}
Suppose that Assumption \ref{ass2.1} holds true. Then without changing their
(own but may be not joint) distributions the processes $\Xi_i(t),\, t\geq 0,\,
i=1,...,\ell$ can be redefined on a richer probability space where there exist 
also standard Brownian motions $B_i(t),\, t\geq 0,\, i=1,...,\ell$ such that 
for some constants $\al>0$ and $\sig_i\geq 0,\, i=1,...,\ell$,
\begin{equation}\label{2.22}
\Xi_i(t)-\sig_iB_i(t)\ll t^{\frac 12-\al}\,\,\mbox{a.s.}.
\end{equation}
\end{theorem}

As usual (see \cite{PS} and \cite{HH}), relying on the well known invariance
principle in the law of iterated logarithm for the Brownian motion (see
\cite{St0}) we obtain immediately from the above theorem the following result.
\begin{corollary}\label{corlil}
Let $K_i$ be the compact set of absolutely continuous functions $x$ in 
$C[0,1]$ with $x(0)=0$ and $\int_0^1\dot x^2(u)du\leq\sig_i^2$ and set
$\zeta_{i,t}(u)=(2t\ln\ln t)^{-1/2}\Xi_i(tu)$, $u\in[0,1]$. Then
the family $\zeta_{i,t},\, t\geq 3$ is relatively compact in the topology
of uniform convergence and as $t\to\infty$ the set of all a.s. limit points
of $\zeta_{i,t}$ coincides with $K_i$. Let $K$ be the compact set of
functions $x\in C[0,1]$ which can be written in the form
 $x(u)=\sum_{1\leq i\leq\ell}x_i(u)$ with
 $x_i\in K_i,\, i=1,...,\ell$. Set $\zeta_t(u)=(2t\ln\ln t)^{-1/2}
 \Xi(tu)$, $u\in[0,1]$. Then 
the family $\zeta_t,\, t\geq 3$ is relatively compact in the topology
of uniform convergence and as $t\to\infty$ the set of all a.s. limit points
of $\zeta_t$ is contained in $K$.
\end{corollary}

In order to understand our assumptions observe that $\varpi_{q,p}$
 is non-increasing in $q$ and non-decreasing in $p$. Hence,
 for any pair $p,q\geq 1$,
 \[
 \varpi_{q,p}(n)\leq\psi(n).
 \]
 Furthermore, by the real version of the Riesz--Thorin interpolation 
 theorem (see, for instance, \cite{Ga}, Section 9.3) if 
 $\theta\in[0,1],\, 1\leq p_0,p_1,q_0,q_1\leq\infty$ and
 \[
 \frac 1p=\frac {1-\theta}{p_0}+\frac \theta{p_1},\,\,\frac 1q=\frac
 {1-\theta}{q_0}+\frac \theta{q_1}
 \]
 then
 \begin{equation*}
\varpi_{q,p}(n)\le 2(\varpi_{q_0,p_0}(n))^{1-\theta}
(\varpi_{q_1,p_1}(n))^\theta.
\end{equation*}
Since, clearly, $\varpi_{q_1,p_1}\leq 2$ for any $q_1\geq p_1$ it follows
for pairs $(\infty,1)$, $(2,2)$ and $(\infty,\infty)$ that for all 
$q\geq p\geq 1$,
 \begin{eqnarray*}
&\varpi_{q,p}(n)\le (2\alpha(n))^{\frac{1}{p}-\frac{1}{q}},\,
 \varpi_{q,p}(n)\le 2^{1+\frac 1p-\frac 1q}(\rho(n))^{1-\frac 1p+\frac 1q}\\
&\mbox{and}\,\,\varpi_{q,p}(n)\le 2^{1+\frac 1p}(\phi(n))^{1-\frac 1p}.
\end{eqnarray*}
We observe also that by the H\" older inequality for $q\geq p\geq 1$
and $\alpha\in(0,p/q)$,
\begin{equation*}
\beta(q,r)\le 2^{1-\alpha}  [\beta(p,r)]^\alpha \gamma^{1-\al}_{\frac{pq(1-\al)}
{p-q\al}}
\end{equation*}
with $\gamma_\theta$ defined in (\ref{2.18}). Thus, we can formulate 
Assumption \ref{ass2.1} in terms of more familiar $\alpha,\,\rho,\,\phi,$
and $\psi$--mixing coefficients and with various moment conditions.

The strategy of the proof of Theorem 2.2 consists of several steps. First, 
we split the sum $\Xi_i(t)$ into a sum of "big" and "small" growing blocks
so that the total contribution of small block can be disregarded and their 
sole purpose is to provide sufficient separation between big blocks. Growing
blocks will enable us to approximate their members by conditional expectations
as in (\ref{2.3}) with increasing precision which differs from \cite{KV} and
is an important point in obtaining our estimates. In spite
of the fact that big blocks still remain strongly dependent in our setup
the technique of \cite{KV} enables us to treat them as if they were weakly
dependent. Namely, employing appropriate estimates from \cite{KV} we construct
a martingale approximation of sums of big blocks with an error sufficient
for our purposes. Finally, we rely on the Skorokhod embedding of martingales
into a Brownian motion and estimate the distance between the embedded 
process and the Brownian motion. 

We observe that though the Skorokhod embedding preserves distribution of each
one dimensional martingale it does not preserve, in general, joint 
distributions of several martingales when we employ it simultaneously to
$\ell$ of them as in our case. By this reason we obtain strong approximations
(\ref{2.22}) for each $\Xi_i(t)$ but we do not obtain a strong approximation
of the sum $\Xi(t)$ by the Gaussian process $\sum_{i=1}^\ell \sig_iB_i(t)$ 
which according to \cite{KV} is the weak limit of processes $N^{-1/2}\Xi(Nt)$
as $N\to\infty$. In fact, this is connected with multidimensional strong
approximation theorems where the Skorokhod embedding is not applicable
while other methods employed usually in these circumstances do not seem to
work for in our nonconventional setup.

\section{blocks and martingale approximation}\label{sec3}
\setcounter{equation}{0}

The following result which is a part of Corollary 3.6 from \cite{KV} (improving
in several respects Lemma 3.1 from \cite{Ki}) will be a base for our estimates.

\begin{proposition}\label{prop3.1}
Let $\cG$ and $\cH$ be $\sig$-subalgebras on a probability space $(\Om,\cF,P)$,
$X$ and $Y$ be $d$-dimensional random vectors and $f=f(x,\om),\, x\in\bbR^d$ be a 
collection of random variables measurable with respect to $\cH$ and satisfying
\begin{equation}\label{3.1}
\|f( x,\omega)-f( y,\omega)\|_{q}\le C_1 (1+|x|^\iota + |y|^\iota)|x-y|^\ka
\,\,\mbox{and}\,\,\|f(x,\omega)\|_{q}\le C_2 (1+|x|^\iota)
\end{equation}
where  $g\geq 1$. Set $g(x)=Ef(x,\om)$. Then
\begin{equation}\label{3.3}
\| E(f(X,\cdot)|\cG)-g(X)\|_\up\leq c(1+\| X\|^{\io+2}_{b(\io+2)})
(\vp_{q,p}(\cG,\cH)+\| X-E(X|\cG)\|^\del_{q}),
\end{equation}
 provided $\frac 1\up\geq\frac 1p+\frac 1{b}+\frac {\del}q$ and
 $\ka-\frac dp>\del>0$ with $c=c(C_1,C_2,\io,\io',\ka,
 \del,p,q,\up,d)>0$ depending only on parameters in brackets. Moreover, let
$x=(v,z)$ and $X=(V,Z)$, where $V$ and $Z$ are $d_1$ and $d-d_1$-dimensional
random vectors, respectively, and let $f(x,\om)=f(v,z,\om)$ satisfy (\ref{3.1})
in $x=(v,z)$. Set $\tilde g(v)=Ef(v,Z(\om),\om)$. Then
\begin{eqnarray}\label{3.4}
&\| E(f(V,Z,\cdot)|\cG)-\tilde g(V)\|_\up\leq c(1+\| X\|^{\io+2}
_{b(\io+2)})\\
&\times\big(\vp_{q,p}(\cG,\cH)
+\| V-E(V|\cG)\|^\del_{q}+\| Z-E(Z|\cH)\|^\del_{q}\big).\nonumber
\end{eqnarray}
\end{proposition}

We will use the following notations
\begin{eqnarray}\label{3.6}
&F_{i,r,n}(x_1,x_2,\ldots, x_{i-1},\omega)
=E\big(F(x_1,x_2,\ldots, x_{i-1},X(n))|\cF_{n-r,n+r}\big),\\
&X_r(n)=E\big(X(n)|\cF_{n-r,n+r}\big),\,\,Y_i(q_i(n))=F_i(X(q_1(n)),\ldots,
 X(q_i(n)))\quad\mbox{and}\nonumber\\
&Y_i(j)=0\quad{\rm if}\quad j\ne q_i(n)\quad\mbox{for any}\,\, n,\,\,
Y_{i,r}( q_i(n))=F_{i,r,q_i(n)} (X_r(q_1(n)),\nonumber\\
&\ldots, X_r(q_{i-1}(n)),\omega)\quad\mbox{and}\quad Y_{i,r}(j)=0 
\quad{if}\quad j\ne q_i(n)\quad\mbox{for any}\,\, n.\nonumber
\end{eqnarray}
Next, we fix some positive numbers $4\eta<2\te<\tau<1/2$ which will be 
specified later on and following \cite{PS} introduce pairs of "big" and 
"small" increasing blocks defining for each $i$ random variables $V_i(j)$ 
and $W_i(j)$ inductively so that
\begin{eqnarray}\label{3.7}
&V_i(1)=Y_{i,1}(q_i(1)),\, W_i(1)=Y_{i,1}(q_i(2)),\, a(1)=0,\, b(1)=1\,\,
\mbox{and for}\,\, j>1,\\
& a(j)=b(j-1)+[(j-1)^\te],\,\, b(j)=a(j)+[j^\tau],\,\, r(j)=[j^\eta],
\nonumber\\
&V_i(j)=\sum_{a(j)< l\leq b(j)}Y_{i,r(j)}(q_i(l))\,\,
\mbox{and}\,\, W_i(j)=\sum_{b(j)<l\leq a(j+1)}Y_{i,r(j)}(q_i(l)).\nonumber
\end{eqnarray}
Observe that unlike \cite{KV} but following \cite{PS} the parameter $r(j)$
grows with $j$ increasing precision of conditional expectations approximations.
Let $\nu_i(t)=\max\{ j:\, b(j)+[j^\te]\leq t+1\}$ which is the number of
full small blocks in the sum $\Xi_i(t)$.
We will see that the small blocks $W_i(j),\, j=1,2,...$ make negligible 
contributions to the sum $\Xi_i$ and can be disregarded while the big blocks
$V_i(j),\, j=1,2,...$ are widely separated which enables us to exploit fully
our mixing assumptions. Observe that unlike the sums appearing in standard
limit theorems these big blocks are strongly (and not weakly) dependent but
as in \cite{KV} we will see by means of Proposition \ref{prop3.1} that only
sufficient separation between $q_i(l)$ for different $l$'s plays the role.

Next, set
\begin{equation}\label{3.8}
R_i(m)=\sum_{j=m+1}^\infty E\big( V_i(j)|\cG_m\big)
\end{equation}
and $M_i(m)=V_i(m)+R_i(m)-R_i(m-1)$ where $\cG_m=\cF_{-\infty, q_i(b(m))
+r(m)}$. 
Observe that if $a(j)<l\leq b(j)$ and $j\geq m+1$ then $X=
\big(X_{r(j)}(q_1(l)),...,X_{r(j)}(q_{i-1}(l))\big)$ is 
$\cF_{-\infty,q_{i-1}(l)+r(j)}$ measurable while 
$f(x,\om)=F_{i,r(j),q_i(l)}(x_1,...,x_{i-1},\om)$ is $\cF_{q_i(l)-r(j),\infty}$
measurable. Hence, by (\ref{3.3}) considered with $\cG=\cF_{-\infty,
\max(q_{i-1}(l)+r(j),q_i(b(m))+r(m))}$ and $\cH=\cF_{q_i(l)-r(j),\infty}$ we
obtain that
\begin{equation}\label{3.9}
\| E\big(Y_{i,r(j)}(q_i(l))|\cG_m)\|_{2+\del}\leq C\vp_{q,p}(d_{i,j}(l))
\end{equation}
where $p$ and $q$ satisfy conditions of Proposition \ref{3.1} with $\up=2+\del$
and Assumption \ref{ass2.1},
$C>0$ does not depend on $i,j,l,m$ and
\begin{eqnarray}\label{3.10}
&d_{i,j}(l)=\min(q_i(l)-q_{i-1}(l)-2r(j),\, q_i(l)-q_i(b(m))-r(j)-r(m))\\
&\geq l-b(m)-2r(j)\geq a(j)-b(m)-2r(j)\nonumber
\end{eqnarray}
taking into account that under our assumptions
\begin{equation}\label{3.11}
q_i(l)-q_{i-1}(l)\geq l\,\,\mbox{and}\,\, q_i(l)-q_i(m)\geq l-m
\end{equation}
provided $l$ is large enough. Thus, for $j\geq m+1$,
\begin{equation}\label{3.12}
\big\| E(V_i(j)|\cG_m)\big\|_{2+\del}\leq C\sum_{a(j)<l\leq b(j)}\vp_{q,p}
(l-b(m)-2r(j))
\end{equation}
and
\begin{equation}\label{3.13}
\| R_i(m)\|_{2+\del}\leq C\sum^\infty_{j=m+1}\sum_{a(j)<l\leq b(j)}
\vp_{q,p}(l-b(m)-2r(j))\leq\tilde C<\infty
\end{equation}
for some constant $\tilde C>0$. In particular, the series (\ref{3.8})
converges in $L^{2+\del}(\Om,\cF,P)$, and so the definition of $R_i(m)$ makes
sense. Observe also that $M_i(m)$ is $\cG_m$ measurable and
\[
E(M_i(m)|\cG_{m-1})=E(V_i(m)+R_i(m)|\cG_{m-1})-R_i(m-1)=0
\]
which means that $(M_i(m),\cG_m),\, m=1,2,...$ is a martingale difference
sequence.

We are going to replace the sum $\Xi_i(t)$ by the martingale $\sum_{1\leq
m\leq\nu_i(t)}M_i(m)$ and it will be crucial for our purposes to estimate
the corresponding error. In order to make the first step in this direction
we set
\[
I_1(m)=\sum_{1\leq j\leq m}(V_i(j)-M_i(j))
\]
and relying on (\ref{3.13}) it follows that
\begin{equation}\label{3.14}
\big\| I_1(m)\big\|_2=\| R_i(\nu_i(m))\|_2+\| R_i(0)\|_2\leq 2\tilde C
\end{equation}
for some constant $\tilde C>0$. By Chebyshev's inequality
\begin{equation}\label{3.15}
P\{ |I_1(m)|\geq \frac 12m^{\frac 12+\ve}\}\leq 16\tilde C^2m^{-(1+2\ve)}.
\end{equation}
Observe that
\[
[t]\geq\sum_{1\leq j\leq\nu_i(t)}[j^\tau]\geq\int_0^{\nu_i(t)}u^\tau du=
(1+\tau)^{-1}(\nu_i(t))^{1+\tau},
\]
and so
\begin{equation}\label{3.16}
\nu_i(t)\leq \big((1+\tau)[t]\big)^{1/1+\tau}\leq 2[t]^{1/1+\tau}.
\end{equation}
Hence, taking $m=\nu_i([t])=\nu_i(t)$ and $\ve=\frac 12(\frac 12\tau-\te)
+\frac 12\tau(\frac 12-\te)>0$ we obtain by (\ref{3.15}) and (\ref{3.16})
that
\begin{eqnarray}\label{3.17}
&P\big\{ \big\vert I_1(\nu_i([t]))\big\vert\geq [t]^{\frac 12(1-\te)}\big\}
\leq P\big\{ |I_1(\nu_i([t]))|\\
&\geq(\frac 12\nu_i([t]))^{\frac 12+\ve}\big\}
\leq 16\tilde C^2(\nu_i([t]))^{-(1+2\ve)}.\nonumber
\end{eqnarray}
Therefore, by the Borel-Cantelli lemma we conclude that for some $n_0=
n_0(\om)$ and all $\nu_i(t)\geq n_0$,
\begin{equation}\label{3.18}
|I_1(\nu_i(t))|\leq t^{\frac 12(1-\te)}\quad\mbox{a.s.}
\end{equation}
Observe that
\[
t<2\sum_{1\leq j\leq\nu_i(t)+1}[j^\tau]\leq\int_0^{\nu_i(t)+1}(u^\tau+1)du
\leq 4\tau^{-1}(\nu_i(t)+1)^{1+\tau},
\]
and so
\begin{equation}\label{3.19}
\nu_i(t)\geq (\tau t/2)^{1/1+\tau}-1.
\end{equation}
Hence, if $t\geq 2\tau^{-1}(n_0+1)^{1+\tau}$ then (\ref{3.18}) holds true.

Next, set
\[
I_2(m)=\sum_{a(m+1)\leq l<a(m+2)}|Y_i(q_i(l))|.
\]
Since $a(\nu_i(t)+1)\leq t<a(\nu_i(t)+2)$ then
\begin{equation}\label{3.20}
\big\vert\sum_{a(\nu_i(t))\leq l\leq t}Y_i(q_i(l))\big\vert\leq 
I_2(\nu_i(t)).
\end{equation}
By (\ref{2.6}) and (\ref{2.21}),
\begin{equation}\label{3.20+}
\| Y_i(q_i(l))\|_{2+\del}\leq C<\infty
\end{equation}
for some $C>0$ independent of $l$ and since
by the construction $a(m+2)-a(m+1)=(m+1)^\tau+(m+1)^\te$ we see that
\begin{equation}\label{3.21}
\| I_2(m)\|_{2+\del}\leq\sum_{a(m+1)\leq l<a(m+2)}\| Y_i(q_i(l))\|_{2+\del}
\leq 2C(m+1)^\tau.
\end{equation}
 By (\ref{3.16}) and Chebyshev's inequality
\begin{eqnarray}\label{3.22}
&P\big\{ |I_2(\nu_i(t))|\geq [t]^{\frac 12(1-\ve)}\big\}\\
&\leq P\big\{|I_2(\nu_i(t))|\geq(\frac 12\nu_i(t))^{\frac 12(1+\tau)
(1-\ve)}\big\}\leq \tilde C(\nu_i(t))^{-1-\be}\nonumber
\end{eqnarray}
for some $\tilde C>0$ independent of $t$ where we assume that 
$\tau\leq\frac 14\min(\del,1)$ and take $\ve=
\frac 1{16}\min(\del,1),\,\be=\frac 1{128}\min(\del,1)$. As in (\ref{3.18})
we conclude using (\ref{3.19}) and the Borel-Cantelli lemma that
\begin{equation}\label{3.23}
I_2(\nu_i(t))\leq t^{\frac 12(1-\ve)}\quad\mbox{a.s.}
\end{equation}
for all $t\geq t_0$ and some random variable $t_0=t_0(\om)<\infty$.

Next, we estimate contribution of the small blocks. Let $l>j$ then
$W_i(j)$ is measurable with respect to $\cG=\cF_{-\infty,q_i(a(j+1))+r(j)}$
provided $j^\eta>2$, and so applying (\ref{3.3}) with such $\cG$,
$f(x_1,...,x_{i-1},\om)=F_{i,r(l),q_i(n)}(x_1,...,x_{i-1},\om)$ where
$b(l)<n\leq a(l+1)$, $\cH=\cF_{q_i(b(l))-r(l),\infty}$ we obtain by 
(\ref{2.6}), (\ref{2.21}) and (\ref{3.11}) that for $n$ large enough,
\begin{eqnarray}\label{3.24}
&\big\vert EW_i(j)Y_{i,r(l)}(q_i(n))\big\vert =\big\vert 
E\big(W_i(j)E(Y_{i,r(l)}(q_i(n))|\cG)\big)\big\vert\\
&\leq C_1\vp_{q,p}(n-r(l)-a(j+1)-r(j))\| W_i(j)\|_2\nonumber
\end{eqnarray}
for some $C_1>0$ independent of $j,n,l$ satisfying the conditions above.
Since by (\ref{2.6}), (\ref{2.21}) and the definition of blocks,
\begin{equation}\label{3.25}
\| W_i(j)\|_2\leq\sum_{b(j)<l\leq a(j+1)}\| Y_{i,r(j)}(q_i(l))\|_2\leq 
C_2[j^\te]
\end{equation}
for some $C_2>0$ independent of $j$, then
\begin{equation}\label{3.26}
|E\big(W_i(j)W_i(l)\big)|\leq C_1C_2[j^\te][l^\te]\vp_{q,p}
(\sum_{j<m\leq l}([m^\tau]-2[m^\eta])).
\end{equation}
Hence, by (\ref{2.19}), (\ref{3.25}) and (\ref{3.25}) for any positive 
integers $m<n$,
\begin{eqnarray}\label{3.27}
&E\big(\sum_{m<l\leq n}W_i(l)\big)^2\leq\sum_{m<l\leq n}\big(W_i^2(l)+
2\sum_{m<j<l}|E(W_i(j)W_i(l))|\big)\\
&\leq C_3\sum_{m<l\leq n}l^{2\te}\leq C_4(n^{1+2\te}-m^{1+2\te})\nonumber
\end{eqnarray}
for some $C_3,C_4>0$ independent of $m$ and $n$. It follows by Theorem A1
from \cite{PS} together with (\ref{3.16}) that
\begin{equation}\label{3.28}
\big\vert\sum_{1\leq m\leq\nu_i(t)}W_i(j)\big\vert\ll (\nu_i(t))^{\frac 12+\te}
\log^3\nu_i(t)\le 2t^{\frac 12-\ve}\quad\mbox{a.s.}
\end{equation}
where $\ve<(\frac 12\tau-\te)(1+\tau)^{-1}$ and $t$ is large enough.

Next, set
\[
I_3(m)=\big\vert\sum_{1\leq j\leq m}\sum_{a(j)<l\leq a(j+1)}\big(Y_i(q_i(l))-
Y_{i,r(j)}(q_i(l))\big)\big\vert
\]
By (\ref{2.3}), (\ref{2.5}), (\ref{2.20}) and H\" older's inequality (see 
Lemma 4.11 in \cite{KV}),
\begin{equation}\label{3.29}
\| Y_i(q_i(l))-Y_{i,r(j)}(q_i(l))\|_2\leq C\be^\del_q(r(l))
\end{equation}
for some $q,\del>0$ satisfying (\ref{2.21}) and for a constant $C>0$
independent of $j$. Hence, by (\ref{2.20}),
\begin{equation}\label{3.30}
\|I_3(\nu_i(t)\|_2\le\tilde C<\infty
\end{equation}
for some constant $\tilde C>0$ independent of $t$. Proceeding in the same 
way as in (\ref{3.18}) we obtain that for some random variable $t_0=t_0(\om)$,
\begin{equation}\label{3.31}
|I_3(\nu_i(t))|\leq t^{\frac 12(1-\te)}\quad\mbox{a.s.}
\end{equation}
whenever $t\geq t_0$. Finally, collecting (\ref{3.18}), (\ref{3.23}),
(\ref{3.28}) and (\ref{3.31}) we conclude that
\begin{equation}\label{3.32}
|\Xi_i(t)-\sum_{1\leq j\leq\nu_i(t)}M_i(j)|\ll t^{\frac 12-\ve}
\end{equation}
for some $\ve>0$.

\section{Completing the proof via Skorokhod embedding}\label{sec4}
\setcounter{equation}{0}

A martingale version of the Skorokhod embedding (representation) theorem
(see \cite{St}, Theorem 4.3 and \cite{HH}, Theorem A1) applied to our
martingale $\cM_i(m)=\sum_{1\leq j\leq m}M_i(j)$ yields that if $\{ B_i(t),
\, t\geq 0\}$ is a standard Brownian motion then there exist non-negative
random variables $T_j=T_{i,j}$ such that the processes 
\begin{equation}\label{4.1}
\{ B_i(\sum_{1\leq j\leq m}T_j),\, m\geq 1\}\quad\mbox\quad\{\cM_i(m),\,
m\geq 1\}
\end{equation}
have the same distributions. Hence, without loss of generality we can 
redefine $\{ M_i(j),\, j\geq 1\}$ by
\begin{equation}\label{4.2}
M_i(m)=B_i(\sum_{1\leq j\leq m}T_j)-B_i(\sum_{1\leq j\leq m-1}T_j)
\end{equation}
and can keep the same notations for both $M_i(m)$ and $\cM_i(m)$. In fact,
we will redefine also the processes $X(n),\, V_i(m),\, W_i(m)$ we had before
on a richer and common with $M_i(m)$ probability space so that all marginal
and joint distributions remain intact. Furthermore, the embedding theorem 
cited above yields that if $\cA_m$ is the $\sig$-algebra generated by
$\{ B_i(t),\, 0\leq t\leq\sum_{1\leq j\leq m}T_j\}$ then $T_m$ is $\cA_m$
measurable, $B_i(\sum_{1\leq j\leq m}T_j+s)-B_i(\sum_{1\leq j\leq m}T_j)$
is independent of $\cA_m$ for any $s>0$,
\begin{equation}\label{4.3}
E(T_m|\cA_{m-1})=E(M^2_i(m)|\cA_{m-1})=E(M^2_i(m)|\cG_{m-1})=
E(M^2_i(m)|\tilde\cG_{m-1})
\end{equation}
and
\begin{equation}\label{4.4}
E(T^u_m|\cA_{m-1})\leq c_uE(|M_i(m)|^{2u}|\cA_{m-1})
\end{equation}
where $c_u>0$ depends only on $u\geq 1$, $\cA_m\supset\cG_m\supset\tilde\cG_m
=\sig\{ M_i(j),\, 1\leq j\leq m\}$ and $\cG_m$ is the same as in (\ref{3.8}).

In order to exploit the representation
\begin{equation}\label{4.5}
\cM_i(m)=B_i(\sum_{1\leq j\leq m}T_j)
\end{equation}
we have to establish a strong law of large numbers with appropriate error
estimates for sums of $T_j$'s in the form
\begin{equation}\label{4.6}
\big\vert\sum_{1\leq t\leq\nu_i(t)}T_j-\sig_i^2t\big\vert=O(t^{1-\la})\quad
\mbox{a.s.}
\end{equation}
for some $\la>0$ and $\sig_i\geq 0$. This would imply that
\begin{equation}\label{4.7}
\big\vert B_i(\sum_{1\leq j\leq\nu_i(t)}T_j)-B_i(\sig_i^2t)\big\vert\ll 
t^{\frac 12-\tilde\la}\quad\mbox{a.s.}
\end{equation}
for some $\tilde\la<\frac 12\la$. Indeed, set $\tau_i(t)=\sum_{1\leq j
\leq\nu_i(t)}T_j$. Then (\ref{4.6}) means that $|\tau_i(t)-\sig_i^2t|\leq
Qt^{1-\la}$ for some random variable $Q=Q(\om)<\infty$ a.s. Introducing
the events $\Om_N=\{ Q\leq N\}$ we obtain
\[
A(t)=|B_i(\tau_i(t))-B_i(\sig_i^2t)|\bbI_{\Om_N}\leq A_1(t)+A_2(t)
+A_3(t)
\]
where
\begin{eqnarray*}
&A_1(t)=\sup_{0\leq s\leq Nt^{1-\la}}|B_i(\sig_i^2t+s)-B_i(\sig_i^2t)|,\\
&A_2(t)=|B_i(\sig_i^2t)-B_i(\sig_i^2t-Nt^{1-\la})|\,\,\,\mbox{and}\\
& A_3(t)=\sup_{0\leq s\leq Nt^{1-\la}}|B_i(\sig_i^2t
-Nt^{1-\la}+s)-B_i(\sig_i^2t-Nt^{1-\la})|.
\end{eqnarray*}
By the martingale moment inequalities for the Brownian motion
\[
EA^{2m}_j(t)\leq C_mN^mt^{m(1-\la)},\,\, j=1,2,3
\]
where $C_m>0$ depends only on $m\geq 1$. Thus
\[
P\{ A(n)>n^{\frac 12-\tilde\la}\}\leq 3^{2m-1}C_mN^mn^{-m(\la-2\tilde\la)}.
\]
Choose $\tilde\la<\frac 12\la$ and $m\geq 2(\la-2\tilde\la)^{-1}$ then
$n^{-m(\la-2\tilde\la)}\leq n^{-2}$, and so the probabilities above form
a converging series. Hence, by the Borel--Cantelli lemma there exists 
$n_0=n_0(\om)<\infty$ such that
\[
A(n)\leq n^{\frac 12-\tilde\la}\,\,\,\mbox{ a.s. for all }\,\, n\geq n_0.
\]
Since $\nu_i(t)$, and so also $\tau_i(t)$, can change only at integer $t$
and since $\Om_N\uparrow\tilde\Om$ as $N\uparrow\infty$ with $P(\tilde\Om)
=1$ we conclude that, indeed, (\ref{4.6}) implies (\ref{4.7}).
Finally, redefining without changing distributions
all processes once again we can replace $B_i(\sig^2t)$ by $\sig_iB_i(t)$
arriving at the assertion of Theorem \ref{thm2.2}.

We start deriving (\ref{4.6}) by writing
\begin{equation}\label{4.8}
\sum_{1\leq j\leq m}(T_j-M^2_i(j))=D^{(1)}(m)-D^{(2)}(m),
\end{equation}
where
\[
D^{(1)}(m)=\sum_{1\leq j\leq m}(T_j-E(T_j|\cA_{j-1})),\quad D^{(2)}(m)=
\sum_{1\leq j\leq m}(M^2_i(j)-E(M^2_i(j)|\cG_{j-1})),
\]
and using (\ref{4.3}) in order to have (\ref{4.8}). Set $R^{(1)}(j)=T_j-
E(T_j|\cA_{j-1})$ then $(R^{(1)}(j),\,\cA_j)_{j\geq 1}$ is a martingale 
differences sequence. By (\ref{3.13}), (\ref{3.20+}) and (\ref{4.4}) for 
any $j\geq 1$,
\[
E|R^{(1)}(j)|^{1+\frac 12\del}\leq 2E|T_j|^{1+\frac 12\del}\leq 
2c_{1+\frac 12\del}E|M_i(j)|^{2+\del}\leq C(1+E|V_i(j)|^{2+\del})\leq
\tilde Cj^{(2+\del)\tau}
\]
for some $C,\tilde C>0$ independent of $j$. Observe that 
$(j^{-(1+\tau+\ve)}R^{(1)}(j),\cA_j)_{j\geq 1}$ is also a martingale
 differences sequence and assume that $\tau\leq\del/4$ and $\tau+\ve\leq 1/4$.
 Then
\[
\sum_{j=1}^\infty j^{-(1+\frac 12\del)(1+\tau-\ve)}E|R^{(1)}(j)|^{1+
\frac 12\del}\leq\tilde C\sum_{j=1}^\infty j^{-(1+\frac \del 8)}<\infty,
\]
and so by the standard result on martingale series (see Theorem 2.17   
in \cite{HH}), 
\[
\sum_{1\leq j\leq\infty}j^{-(1+\tau-\ve)}R^{(1)}(j)\quad\mbox{converges a.s.}
\] 
Hence, by Kronecker's lemma 
\[
m^{-(1+\tau-\ve)}\sum_{j=1}^mR^{(1)}(j)=m^{-(1+\tau-\ve)}D^{(1)}(m)\to 0
\quad\mbox{a.s.}\quad\mbox{as}\,\,\, m\to\infty,
\]
and so by (\ref{3.16}),
\begin{equation}\label{4.9}
t^{-(1-\frac {\ve}{1+\tau})}|D^{(1)}(\nu_i(t))|\leq 4(\nu(t))^{-(1+\tau-\ve)}
|D^{(1)}(\nu_i(t))|\to 0\,\,\,\mbox{a.s.}\,\,\,\mbox{as}\,\, t\to\infty.
\end{equation}
Setting $R^{(2)}(j)=M^2_i(j)-E(M^2_i(j)|\cG_{j-1})$ we obtain that 
$(R^{(2)}(j),\,\cG_j)_{j\geq 1}$ is a martingale differences sequence, as well,
 and by (\ref{3.13}) and (\ref{3.20+}),
\[
E|R^{(2)}(j)|^{1+\frac 12\del}\leq 2E|M_i(j)|^{2+\del}\leq 
C(1+E|V_i(j)|^{2+\del})\leq\tilde Cj^{(2+\del)\tau}
\]
for some $C,\tilde C>0$ independent of $j$. Thus, in the same way as above, we 
see that
\begin{equation}\label{4.10}
t^{-(1-\frac {\ve}{1+\tau})}|D^{(2)}(\nu_i(t))|\to 0\,\,\,\mbox{a.s.}\,\,\,
\mbox{as}\,\, t\to\infty.
\end{equation}

It follows from (\ref{4.8})--(\ref{4.10}) that in order to obtain (\ref{4.6})
it suffices to show that there exists $\sig_i\geq 0$ such that
\begin{equation}\label{4.11}
\big\vert\sum_{1\leq j\leq\nu_i(t)}M_i^2(j)-\sig_i^2t\big\vert=O(t^{1-\la})
\quad\mbox{a.s.}
\end{equation}
for some $\la>0$. By the definition of $M_i(j)$ and the Cauchy inequality,
\begin{equation}\label{4.12}
\big\vert\sum_{1\leq j\leq m}(M^2_i(j)-V_i^2(j))\big\vert\leq(A_i(m))^{1/2}
\big(2(\sum_{1\leq j\leq m}V_i^2(j))^{1/2}+(A_i(m))^{1/2}\big)
\end{equation}
where $A_i(m)=\sum_{1\leq j\leq m}\rho_i^2(j)$, $\rho_i(j)=R_i(j)-R_i(j-1)$
and $E|A_i(m)|\leq m\tilde C^2$ by (\ref{3.13}). Fix $\be>0$ and for each
$l\geq 1$ set $m_l=[l^{2/\be}]$ then by Chebyshev's inequality
\begin{equation}\label{4.13}
P\{ |A_i(m)|\geq m_l^{1+\be}\}\leq\tilde C^2m_l^{-\be}\leq\tilde{\tilde C}
l^{-2}
\end{equation}
for some $\tilde{\tilde C}>0$ independent of $l$. Therefore, by the
Borel-Cantelli lemma for all $l\geq l_0=l_0(\om)<\infty$,
\begin{equation*}
|A_i(m_l)|<m^{1+\be}_l\quad\mbox{a.s.}
\end{equation*}
If $m_l\leq\nu_i(t)<m_{l+1}$ and $l\geq l_0$ then by (\ref{3.16}),
\[
|A_i(\nu_i(t))|\leq |A_i(m_{l+1})|<m^{1+\be}_{l+1}<(\nu_i(t))^{1+\be}
(\frac {m_{l+1}}{m_l})^{1+\be}\leq Ct^{\frac {1+\be}{1+\tau}}\,\,\mbox{a.s.}
\]
where $C>0$ does not depend on $l$. Choosing  $\be=\tau/2$ we obtain that
\begin{equation}\label{4.14}
A_i(\nu_i(t))\leq Ct^{1-\frac {\tau}{2(1+\tau)}}\quad\mbox{a.s.}
\end{equation}
for all $t\geq t_0=t_0(\om)<\infty$ where in view of (\ref{3.19}) we can
take $t_0=\frac 2\tau((l_0+1)^{2/\be}+2)^{1+\tau}$.

It follows from (\ref{4.12}) and (\ref{4.14}) that in order to obtain
(\ref{4.11}) it remains to show that 
\begin{equation}\label{4.15}
\big\vert\sum_{1\leq j\leq\nu_i(t)}V_i^2(j)-\sig_i^2t\big\vert =O(t^{1-\la})
\quad\mbox{a.s.}
\end{equation}
for some $\la>0$. Next, we will make yet another reduction showing that
(\ref{4.15}) will follow if 
\begin{equation}\label{4.16}
\big\vert\big(\sum_{1\leq j\leq t}Y_i(q_i(l))\big)^2-\sig_i^2t\big\vert
=O(t^{1-\la})\quad\mbox{a.s.}
\end{equation}
for some $\la>0$. A transition from (\ref{4.16}) to (\ref{4.15}) proceeds
in the same way as in Lemma 7.3.5 of \cite{PS} but for readers' convenience
we sketch also here the corresponding argument.

First, we write
\begin{equation}\label{4.17}
\big\vert E\big(\sum_{1\leq l\leq t}Y_i(q_i(l))\big)^2-\sum_{1\leq j\leq
\nu_i(t)}V^2_i(j)\big\vert\leq J_1(t)+J_2(\nu_i(t))
\end{equation}
where
\[
J_1(t)=\big\vert E\big(\sum_{1\leq l\leq t}Y_i(q_i(l))\big)^2-
\sum_{1\leq j\leq\nu_i(t)}EV^2_i(j)\big\vert
\]
and
\[
J_2(m)=\big\vert\sum_{1\leq j\leq m}\big(V^2_i(j)-EV_i^2(j)\big)\big\vert.
\]
Next,
\begin{equation}\label{4.18}
J_1(t)\leq J_{11}(\nu_i(t))+J_{12}(\nu_i(t))+J_{13}(\nu_i(t))+J_{14}(\nu_i(t))
+J_{15}(\nu_i(t))
\end{equation}
where
\begin{eqnarray*}
&J_{11}(m)=2\sum_{1\leq j<\tilde j\leq m}|EV_i(j)V_i(\tilde j)|,\,\,\,
J_{12}(m)=E\big(\sum_{1\leq j\leq m}W_i(j)\big)^2,\\
& J_{13}(m)=E(I_2(m))^2,\,\, J_{14}(m)=E(I_3(m))^2\,\,\mbox{and}\\
&J_{15}(m)=2\|\sum_{1\leq j\leq m}V_i(j)\|_2\big(J^{1/2}_{12}(m)+
J^{1/2}_{13}(m)+J^{1/2}_{14}(m)\big)
\end{eqnarray*}
with $I_2$ and $I_3$ the same as in (\ref{3.21}) and (\ref{3.30}), 
respectively. Using (\ref{3.9})--(\ref{3.12}) we obtain similarly to 
(\ref{3.24})--(\ref{3.27}) that
\begin{eqnarray}\label{4.19}
&J_{11}(m)\leq 2C\sum_{1\leq j<\tilde j\leq m}\| V_i(j)\|_2\sum_{a\tilde j)
<l\leq b(\tilde j)}\vp_{q,p}(l-b(j)-2r\tilde j))\\
&2C\sum_{1\leq j<\tilde j\leq m}\tilde j^{\tau}\sum_{a(\tilde j)
<l\leq b(\tilde j)}\vp_{q,p}(l-b(j)-2r\tilde j))\leq\tilde Cm\nonumber
\end{eqnarray}
for some $C,\tilde C>0$ independent of $m$. For $J_{12}(m)$, $J_{13}(m)$
and $J_{12}(m)$ we already have appropriate estimates in (\ref{3.27}),
(\ref{3.21}) and (\ref{3.30}), respectively. Employing (\ref{3.3}) from
Proposition \ref{prop3.1} together with Assumption \ref{ass2.1} in order
to estimate $a_{ln}=|EY_i(q_i(l))Y_i(q_i(n))|$ we see (see (\ref{4.31}) 
and (\ref{4.34}) below as well as Lemma 5.1 from \cite{KV}) that 
$\sum_{1\leq l<n\leq t}a_{ln}$ is of order $O(t)$, and so for $m\leq\nu_i(t)$,
\begin{equation}\label{4.20}
E\big(\sum_{1\leq j\leq m}V_i(j)\big)^2\leq\sum_{1\leq l\leq t}Y^2_i(q_i(l))
+2\sum_{1\leq l<n\leq t}a_{ln}\leq Ct
\end{equation}
for some $C>0$ independent of $t$. Combining (\ref{3.16}), (\ref{3.21}),
(\ref{3.27}), (\ref{3.30}) and (\ref{4.18})--(\ref{4.20}) we obtain that
\begin{equation}\label{4.21}
J_{11}(\nu_i(t))\leq\tilde Ct^{1-\ve}
\end{equation}
for some $\tilde C>0$ independent of $t$ where $\ve=(\tau-2\te)/(1+\tau)$.

In order to estimate $J_2(t)$ we set 
\begin{equation*}
   U_i(j)=\left\{\begin{array}{ll}
  V_i^2(j)-EV^2_i(j) &\mbox{if}\,\,\,\, |V_i^2(j)-EV_i^2(j)|\leq j^{1+\sig}\\
  0 &\mbox{otherwise}
  \end{array}\right.
\end{equation*}
where $\sig\in [\tau,\frac \del 4]$ will be further specified later on. 
Observe that 
\begin{eqnarray}\label{4.22}
& P\{ U_i(j)\ne V_i^2(j)-EV^2_i(j)\} =P\{ |V_i^2(j)-EV_i^2(j)|>j^{1+\sig}\}\\
&\leq 2^{1+\frac \del 2}j^{-(1+\sig)(1+\frac \del 2)}\leq 2^{1+\frac \del 2}
j^{-(1+\frac \del 2)(1+\sig-2\tau)}.\nonumber
\end{eqnarray}
Since $\sig\geq 2\tau$ then the power of $j$ in the right hand side of 
(\ref{4.22}) is less than $-1$, and so by the Borel--Cantelli lemma with
probability one the event $\{ U_i(j)\ne V_i^2(j)-EV^2_i(j)\}$ can occur 
only finite number of times. Hence, the asymptotical behaviot as $m\to
\infty$ of $J_2(m)$ and of $J_3(m)=|\sum_{1\leq j\leq m}U_i(j)|$ is the
same (up to a random variable independent of $m$)
and it suffices to estimate the latter. Set $U^*_i(j)=U_i(j)-EU_i(j)$.
Using (\ref{3.12}) we obtain that for $j<j'$, 
\begin{equation}\label{4.23}
|EU^*_i(j)U^*_i(j')|\leq (jj')^{1+\sig}\vp_{q,p}(j^\te+\sum_{j<m\leq j'}
m^\tau).
\end{equation}
Next,
\begin{equation}\label{4.24}
E(U^*_i(j))^2\leq 2j^{(1+\sig)(1-\frac \del 2)}E|U^*_i(j)|^{1+\frac \del 2}
\leq 32j^{(1+\sig)(1-\frac \del 2)}E|V_i(j)|^{2+\del}\leq Cj^{1+2\tau-\ve}
\end{equation}
for some $C>0$ independent of $j$ where $\ve=\frac \del 2-\sig+
\frac {\del\sig}2-\tau\del$ and we choose $\sig$ and $\tau$ so small that
$\ve\geq\del/8$. It follows from (\ref{2.19}), (\ref{4.23}) and (\ref{4.24})
that for some $\tilde C>0$ independent of $n$ and $m$,
\begin{equation}\label{4.25}
E\big(\sum_{j=m+1}^nU^*_i(j)\big)^2\leq\tilde C(n^{2+2\tau-\ve}-
m^{2+2\tau-\ve})
\end{equation}
and applying again Theorem A1 from \cite{PS} we obtain by (\ref{3.16}) and
(\ref{4.25}) that 
\begin{equation}\label{4.26}
\big\vert\sum_{1\leq j\leq\nu_i(t)}U^*_i(j)\big\vert\ll(\nu_i(t))^{1+\tau-
\frac 12\ve}\leq 2t^{1-\frac {\ve}{2(1+\tau)}}\,\,\mbox{a.s.}
\end{equation}
Hence, $J_2(\nu_i(t))\ll t^{1-\frac {\ve}{2(1+\tau)}}$ a.s., as well.

Finally, it remains to establish (\ref{4.16}). In fact, existence of the
limit
\[
\lim_{t\to\infty}t^{-1}E\big(\sum_{1\leq n\leq t}Y_i(q_i(n))\big)^2=\sig_i^2
\]
and its computation is given in Propositions 4.1 and 4.5 from \cite{KV}
and we only have to explain the estimate (\ref{4.16}) which is actually
hidden inside the proof there. If $i\leq k$ then the above limit has the form
(see Proposition 4.1 in \cite{KV}),
\begin{equation}\label{4.27}
\sig_i^2=\sum_{l=-\infty}^\infty a_i(l)\,\,\mbox{with}\,\, a_i(l)=
\int F_i(x_1,...,x_i)F_i(y_1,...,y_i)\prod_{1\leq u\leq i}d\mu_{ul}(x_u,y_u)
\end{equation}
where $\mu_n$ is the same as in (\ref{2.4}) and $d\mu_0(x,y)=\del_{xy}d\mu(x)$
is the measure supported by the diagonal. If $i> k$ then 
(see Proposition 4.5 in \cite{KV}),
\begin{equation}\label{4.28}
\sig_i^2=\int F^2_i(x_1,...,x_i)d\mu(x_1)\cdots d\mu(x_i).
\end{equation}

We have
\begin{equation}\label{4.29}
E\big(\sum_{1\leq n\leq t}Y_i(q_i(n))\big)^2=\sum_{1\leq n,n'\leq t}b_i(n,n')
\end{equation}
where 
\[
b_i(n,n')=EF_i\big(X(q_1(n)),...,X(q_i(n))\big)
F_i(X(q_1(n')),...,X(q_i(n'))\big).
\]
If $i\leq k$ then for each integer $m$ we consider $b_i(n,n')$ with $in-in'
=im$. Assume that $|im|\leq\frac 14\max(n,n')$ then we can apply (\ref{3.4})
of Proposition \ref{prop3.1} with $\cG=\cF_{-\infty,(i-3/4)\max(n,n')}$, 
$\cH=\cF_{(i-1/2)\max(n,n'),\infty}$, $V=\big(X(n),...,X(i-1)n;X(n'),...,
X((i-1)n')\big)$ and $Z=(X(in),X(in'))$ which gives that
\begin{eqnarray*}
&\big\vert b_i(n,n')-\int EF_i(X(n),...,X((i-1)n),x)F_i(X(n'),\\
&...,X((i-1)n'),y)d\mu_{im}(x,y)\big\vert
\leq C_1\big(\vp_{q,p}(\frac 14\max(n,n'))+\be_q(\frac 14\max(n,n'))\big)
\end{eqnarray*}
for some $C_1>0$ independent of $n,n'$. Repeating these estimates $i$
times we obtain that
\begin{equation*}
|b_i(n,n')-a_i(m)|\leq C_1(\vp_{q,p}(\frac 14\max(n,n'))+
\be_q(\frac 14\max(n,n')).
\end{equation*}

If $|im|>\max(n,n')$, say $im>\max(n,n')$, then applying (\ref{3.4}) of
Proposition \ref{prop3.1} with $\cG=\cF_{-\infty,\max(in',(i-1)n)+
\frac 1{16}n}$, $\cH=\cF_{(i-1/16)n,\infty}$, $V=(X(n),...,X((i-1)n);
X(n'),...,X(in'))$ and $Z=X(in)$ which yields that
\[
|b_i(n,n')|\leq C_2(\vp_{q,p}(n/16)+\be_q(n/16))
\]
for some $C_2>0$ independent of $n$. The same estimate
holds true if $im<-n/4$ with $n'$ in place of $n$, and so we can replace 
above $n$ by $\max(n,n')$. Next, we want to show that a similar estimate
holds true for $a_i(m)$ when $|im|>\max(n,n')$, assuming first that
$im>\max(n,n')$. Since for $in-in'=im$,
\[
a_i(m)=\int EF_i(x_1,...,x_{i-1}, X(in))F_i(y_1,...,y_{i-1},X(in'))
\prod_{1\leq u\leq i-1}d\mu_{um}(x_u,y_u)
\]
we can apply (\ref{3.4}) of Proposition \ref{prop3.1} with $\cG=
\cF_{-\infty,in'+\frac 1{16}n}$, $\cH=\cF_{(i-\frac 1{16})n,\infty}$,
$V=(x_1,...,x_{i-1};y_1,...,y_{i-1},X(in'))$ and $Z=X(in)$ which 
yields that
\[
|a_i(m)|\leq C_3(\vp_{q,p}(n/16)+\be_q(n/16))
\]
for some $C_3>0$ independent of $n$.
If $im<-\frac n4$ then we obtain a similar estimate with $n$ replaced 
by $n'$, and so we can replace $n$ in the above estimate by $\max(n,n')$.
Collecting the above estimates we obtain that if $i\leq k$ and $in-in'
=im$ for an integer $m$ then
\begin{equation}\label{4.30}
|b_i(n,n')-a_i(m)|\leq C_4(\vp_{q,p}(\frac 1{16}\max(n,n'))+
\be_q(\frac 1{16}\max(n,n'))
\end{equation}
for some $C_4>0$ independent of $n$. By (\ref{2.19}), (\ref{2.20}) and
(\ref{4.30}) we obtain that for $i\leq k$,
\begin{equation}\label{4.31}
\big\vert\sum_{1\leq n,n'\leq t}b_i(n,n')-\sig_i^2\big\vert\leq C_5t^{1-\del}
\end{equation}
for some $C_5>0$ independent of $t$.

Next, we consider the case $i\geq k+1$. It follows from (\ref{2.8}) that if
$n\ne n'$ and $\max(n,n')$ is large enough then $|q_i(n)-q_i(n')|\geq
(\max(n,n'))^\del$. Hence, relying on (\ref{3.4}) in Proposition 
\ref{prop3.1} it is easy to see similarly to above that in this case
 \begin{equation}\label{4.32}
 |b_i(n,n')|\leq C_6\big(\vp_{q,p}(\frac 14|q_i(n)-q_i(n')|)+\be_q
 (\frac 14|q_i(n)-q_i(n')|)\big)
 \end{equation}
 for some $C_6>0$ independent of $n$ and $n'$. In order to estimate the
  difference between $b_i(n,n)$ and $\sig_i^2$ from (\ref{4.28}) we use
  that $q_i(n)-q_{i-1}(n)\geq n$ for large $n$ which yields that
  \[
  |b_i(n,n)-\int EF^2_i(X(n),...,X((i-1)n),x)d\mu(x)|\leq C_7
  \big(\vp_{q,p}(\frac 14n)+\be_q(\frac 14n)\big)
  \]
  for some $C_7>0$ independent of $n$ where we rely on 
  (\ref{3.4}) from Proposition \ref{prop3.1} with $\cG=\cF_{-\infty,
  (i-\frac 34)n}$, $\cH=\cF_{(i-\frac 14)n,\infty}$, $V=(X(n),...,
  X((i-1)n))$ and $Z=X(in)$. Repeating this estimate $i$ times we
  obtain that
  \begin{equation}\label{4.33}
  |\sum_{0\leq n\leq t}b_i(n,n)-t\sig_i^2|\leq C_8\sum_{0\leq n\leq t}
  \big(\vp_{q,p}(\frac 14n)+\be_q(\frac 14n)\big)
  \end{equation}
  for some $C_8>0$ independent of $t$. This together with (\ref{2.8}),
  (\ref{2.19}), (\ref{2.20}) and (\ref{4.32}) yields that
  \begin{equation}\label{4.34}
  |\sum_{0\leq n,n'\leq t}b_i(n,n')-t\sig_i^2|\leq C_9t^{1-\del}
  \end{equation}
  for some $C_9>0$ independent of $t$. Finally, (\ref{4.29}), (\ref{4.31})
  and (\ref{4.34}) yields (\ref{4.16}) completing the proof of Theorem 
  \ref{thm2.2}. \qed

\bibliography{matz_nonarticles,matz_articles}
\bibliographystyle{alpha}

\end{document}